\pdfoutput=1
\documentclass{amsart}

\usepackage{amsmath }
\usepackage{bbm,dsfont}
\usepackage{mathrsfs}
\usepackage[all]{xy}
\usepackage{amsfonts}
\usepackage{ifthen}
\usepackage{rotating}
\usepackage{amssymb}
 
\usepackage{tikz}

\usepackage[pdftex]{hyperref}

\newtheorem{definition}{Definition}

\theoremstyle{plain}
\newtheorem{lemma}{Lemma}[section]
\newtheorem{teo}[lemma]{Theorem}
\newtheorem{propo}[lemma]{Proposition}

\theoremstyle{definition}
\newtheorem{example}[lemma]{Example}
\newtheorem{examples}[lemma]{Examples}

\newtheorem{claim*}{Claim}

\theoremstyle{remark}
\newtheorem{remark}[lemma]{Remark}

\newcommand{\Cox}{\operatorname{Cox}}
\newcommand{\Cl}{\operatorname{Cl}}
\newcommand{\Pic}{\operatorname{Pic}}

\newcommand{\Eff}{\operatorname{Eff}}

\newcommand{\Tor}{\operatorname{Tor}}
\newcommand{\Sp}{\operatorname{Spec}}
\newcommand{\Proj}{\operatorname{Proj}}

\newcommand{\rk}{\operatorname{rk}}

\newcommand{\SL}{\operatorname{SL}}
\newcommand{\Osh}{{\mathcal O}}

\newcommand{\cG}{\mathcal{G}}
\newcommand{\cC}{\mathcal{C}}
\newcommand{\cS}{\mathcal{S}}
\newcommand{\cM}{\mathcal{M}}

\newcommand{\irr}{J}

\newcommand{\ff}{\mathbb{F}}
\newcommand{\pp}{\mathbb{P}}
\newcommand{\qq}{\mathbb{Q}}

\newcommand{\cc}{\mathbb{C}}
\newcommand{\zz}{\mathbb{Z}}

\newcommand{\map}{\rightarrow}

\newcommand{\rmap}{\dashrightarrow}

\newcommand{\defi}[1]{\textsf{#1}}

\begin{document}

\title{\bf A survey on Cox rings}
\author{Antonio Laface}
\address{
Departamento de Matem\'atica \newline
Universidad de Concepci\'on \newline 
Casilla 160-C \newline
Concepci\'on, Chile}
\email{alaface@udec.cl}

\author{Mauricio Velasco}
\address{
Department of Mathematics \newline 
1063 Evans Hall \newline
University of California \newline
Berkeley, CA}
\email{velasco@math.berkeley.edu}

\begin{abstract}
We survey the construction of the Cox ring of an algebraic variety $X$ and study the birational geometry of $X$ when its Cox ring is finitely generated.
\end{abstract}
\maketitle

\section*{Introduction}

When studying the geometry of a projective variety $X$ it is often useful to consider its various projective embeddings. In exchange for this flexibility we are left without a unique choice for what ``The coordinate ring" of a projective variety should be. Hu and Keel proposed a candidate in \cite{hk}, which in number theory~\cite{cts} was also known as the {\em universal torsor}: the {\bf Cox ring} of an algebraic variety. This ring generalizes a construction due to Cox in the case of toric varieties~\cite{COX} and can be loosely described (see Section~\hyperref[sec: Coxrings]{\ref{sec: Coxrings}} for a precise definition) as 
$$\Cox(X) := \bigoplus_{L\in\Pic(X)} H^0(X,L)$$
In this note we study some basic examples of {\bf Mori Dream Spaces}, those varieties whose Cox rings are finitely generated. We discuss the restrictions imposed on their geometry by the condition of  finite generation and describe how Cox rings can be used to study the birational geometry of $X$. Our main sources are~\cite{hk} and~\cite{COX}.

The Cox ring of a toric variety is a multigraded polynomial ring (see {
Section~\hyperref[sec:toric]{\ref{sec:toric}}) so all toric varieties are Mori Dream Spaces. Toric varieties play a fundamental role in the theory because, whenever $\Cox(X)$ is finitely generated it admits a presentation $R\map\Cox(X)\map 0$, where $R$ is the Cox ring of some toric variety $Y$. It turns out that $X$ naturally embeds in $Y$ and their birational geometries are closely related. This allows for a combinatorial description of many fundamental invariants of $X$. For example one can describe its nef cone $\overline{NE}^1(X)$ and the decomposition of its effective cone $\Eff(X)$ into {\em Mori chambers}. We will also discuss how $X$ can be obtained as a quotient of an open subset of ${\rm Spec}({\rm Cox(X)})$ by a torus, much like any toric variety is a quotient of an affine space by a torus. In particular this shows that the points of (an open subset of) ${\rm Spec}(\Cox(X))$ can be thought of as homogeneous coordinates on $X$.

This survey is organized as follows. In section 1 we review the construction of the Cox ring of a toric variety $X$ and show that $X$ is a categorical quotient of an affine space by a torus.
We generalize these concepts in Section 2 where we introduce the Mori chambers of $X$.
We define Mori Dream spaces in section 3, while in section 4 we describe some examples of such varieties. Finally in section 5 we survey a geometric method for finding the ideal of relations of a Cox ring once its generators are known.

There are many interesting results in this area and our choice of focusing on the birational aspects has forced us to leave out many important facts such as their factoriality~\cite{bh,EKW} and their remarkable applications for the study of rational points~\cite{Peyre}.  We do not claim originality for the statements written in this survey. Our main sources for the material presented are~\cite{hk},~\cite{COX} and~\cite{LV}. 

\subsection*{Acknowledgments}
It is a pleasure to thank J\"urgen Hausen, Diane MacLagan and Hal Schenck for comments and suggestions which helped us to improve the manuscript. The first author is partially supported by Proyecto FONDECYT Regular 2008, N. proyecto 1080403. The second author is partially supported by NSF grant ${\rm DMS}-0802851$.

\subsection*{Basic notation}

Throughout this paper $k$ is an algebraically closed field.

\section {Toric Varieties and their Cox rings}

An $n$-dimensional toric variety is a normal variety with an open subvariety $T$ isomorphic to the algebraic torus $(k^*)^n$ and such that the action of the torus on itself can be extended to a regular action on $X$ (see~\cite{FULT} or ~\{co} for background on toric varieties).
We describe the Cox ring of a toric variety and show how to reconstruct $X$ from $\Cox(X)$.

\subsection{Notation and definitions}

Let $V$ be a real vector space. A polyhedral cone $\sigma$ in $V$ is the positive span of a finite set of vectors in $V$. A hyperplane $H$ in $V$ is a supporting hyperplane of $\sigma$ if $\sigma$ is contained in one of the two closed half spaces determined by $H$ and $\sigma \cap H$ is non empty. A face $\tau$ of $\sigma$ is the intersection of $\sigma$ with any supporting hyperplane. We use $\tau\preceq\sigma$ to denote that $\tau$ is a face of $\sigma$. The dual cone $\sigma^{\vee}$ is the subset of $V^{*}$ consisting of vectors $w^{*}$ which map $\sigma$ to the nonnegative real numbers. It can be shown that $\sigma^{\vee}$ is a polyhedral cone in $V^{*}$.

Let $N\simeq \mathbb{Z}^r$ be a lattice and let $M=Hom(N,\mathbb{Z})$ be its dual lattice. Denote by $N_{\mathbb{R}}:=N\otimes_\mathbb{Z}\mathbb{R}$ its associated real vector space. A rational polyhedral cone $\sigma$ in $N$ is a cone in $N_{\mathbb{R}}$ generated by vectors in $N$. In this case the dual cone $\sigma^{\vee}$ is also a rational polyhedral cone in $M_{\mathbb{R}}:=M\otimes_{\mathbb{Z}}\mathbb{R}$.

A fan $\Delta$ in $N$ is a collection of rational polyhedral cones closed under $\preceq$ and such that every two cones in $\Delta$ intersect along a common face. Moreover we always assume that the fans are non-degenerate, that is, not contained in any proper subspace of $N_{\mathbb{R}}$. We will now discuss the relationship between toric varieties and fans.

\subsection{Toric varieties and fans}
Let $N\simeq \mathbb{Z}^r$ be a lattice and let $M=Hom(N,\mathbb{Z})$ be its dual lattice. A rational polyhedral fan $\Delta$ in $N$ specifies a toric variety $X(\Delta)$ as follows:
\begin{enumerate}\itemsep 2mm
\item{For each cone $\sigma\in\Delta$ let $U_{\sigma}:=\Sp(k[\sigma^{\vee}\cap M])$.}
\item{For $\sigma\in \Delta$ and each face $\tau\preceq \sigma$ define $(U_{\sigma})_\tau$ to be the spectrum of the localization $k[\sigma^{\vee}\cap M]_{w}$ where $w\in\sigma^{\vee}\cap M$ determines a supporting hyperplane for $\tau$. Note that there is a canonical isomorphism $(U_{\sigma})_\tau\rightarrow U_{\tau}$.}
\item{ For $\sigma_1,\sigma_2$ in $\Delta$ define $(U_{\sigma_1})_{\sigma_2}:=(U_{\sigma_1})_{\sigma_1\cap \sigma_2}$ and let $X(\Delta)$ be the scheme obtained by glueing the schemes $U_{\sigma}$ along the isomorphisms obtained by composing $(U_{\sigma_1})_{\sigma_2}\rightarrow U_{\sigma_1\cap \sigma_2}\rightarrow (U_{\sigma_2})_{\sigma_1}$.}
\end{enumerate}
It follows immediately from the construction that $X(\Delta)$ is separated, integral and normal. Moreover $U_{0}=\Sp(k[M])=(k^{*})^n$ is an open subvariety of $X(\Delta)$ henceforth denoted $T_N$. The action of the torus $T_N$ on itself extends to a regular action on $X$ which is given locally by the map $\phi: k[\sigma^{\vee}\cap M]\rightarrow k[\sigma^{\vee}\cap M]\otimes_k k[M]$ which sends $m$ to $m\otimes m$. Conversely, every toric variety can be specified by a fan as above.

The lattices $M$ and $N$ and the fan $\Delta$ in $N$ have natural geometric interpretations:
\begin{itemize}\itemsep 2mm
\item{The (closed) points of the torus are $k$-algebra homomorphisms $\phi: k[M]\rightarrow k$ and each element $m$ of $M$ corresponds to coordinates on the torus via $\chi^{m}(\phi)=\phi(m)$. These functions are group homomorphisms to $k^*$ and thus $M$ is the character lattice of the torus.}
\item{The elements $n\in N$ correspond to one-parameter subgroups $\lambda_n$. The coordinates of the point $\lambda_n(t)$ are given by $\chi^{m}(\lambda_n(t))=t^{m(n)}$.}
\item{Let $\Delta(1)$ denote the set of rays of  the fan $\Delta$. A ray $v\in \Delta(1)$ determines an irreducible, torus-invariant, codimension one subvariety $Y(v)$ of $X(\Delta)$. More explicitly, if $\sigma$ is any cone containing $v$, the ideal defining $Y(v)$ in $U_{\sigma}$ is generated by the $x^{u}$ in $\sigma^{\vee}\cap M$ for which $\langle u,v\rangle>0$. 
Moreover, the subvarieties $Y(v)$ generate $\Cl(X)$ (the group of Weil divisors modulo principal divisors) and we have the following exact sequence:
\[0\rightarrow M\rightarrow \mathbb{Z}^{\Delta(1)}\rightarrow {\rm Cl}(X)\rightarrow 0\]
 where $\mathbb{Z}^{\Delta(1)}$ is the free group generated by the symbols $Y(v)$ for $v\in \Delta(1)$ and the first morphism takes $u\in M$ to its divisor 
\[
div(x^{u})=\sum_{v\in\Delta(1)}\langle u,v\rangle~ Y(v).
\]
Note that the first map is injective since we are always assuming that our fans are non-degenerate.}
\end{itemize}

\begin{example} Let $X=\mathbb{P}^1\times \mathbb{P}^1$ with torus $(k^*)^2$ embedded by $(\alpha,\beta)\rightarrow ([\alpha:1],[\beta:1])$.
In this case $N\cong \mathbb{Z}^2$ and the action of $T_N$ on $X$ is given by 
\[
(\alpha,\beta)\cdot ([u:v],[s:t])=([\alpha u:v],[\beta s:t]).
\]

The fan $\Delta_0$ of $\mathbb{P}^1\times \mathbb{P}^1$ has four maximal cones with rays $e_1,e_2,-e_1,-e_2\in N$.
Moreover ${\rm Cl}(X)\cong\mathbb{Z}^2$ since we have the exact sequence
\[
\xymatrix@1@R1pt@C20pt{
0\ar[r] & \mathbb{Z}^{2}\ar[rrr]^-{\left(\begin{array}{cc}
$1$ & $0$\\
$0$ & $1$\\
$-1$ & $0$\\
$0$ & $-1$\\
\end{array}\right)} & & &  \mathbb{Z}^4\ar[r] & {\rm Cl}(X)\ar[r] & 0}.
\]
\end{example}

\subsection{Motivation: Toric varieties as quotients}\label{sec:toric}
Perhaps the simplest construction of $\mathbb{P}^n$ is as the quotient of $\mathbb{A}^{n+1}\setminus \{0\}$ by the action of the torus $k^*$. In this section we will survey the article~\cite{COX} where it is shown that a similar construction can be made for every toric variety $X$. This will serve as motivation for the general definition of the Cox ring of an algebraic variety. 

Let $\Delta$ be a fan in $N$ and let $\Delta(1)$ be the set of its one-dimensional cones.

\begin{definition} The Cox ring of a toric variety $X$ is the polynomial ring 
\[\Cox(X):=k[\{v\ |\ v\in \Delta(1)\}]\]
multigraded by $deg(v_i)=[v_i]$, where $[v_i]$ denotes the class of $Y(v_i)$ in ${\rm Cl(X)}$. 
\end{definition}
\begin{definition} For each cone $\sigma\in \Delta$ let $v^{\overline\sigma}$ be the product of the variables $v$ corresponding to rays not contained in $\sigma$. The irrelevant ideal of $Cox(X)$ is the ideal 
\[
\irr_X=(\{x^{\overline\sigma}\ |\ \sigma\in \Delta\}).
\]
\end{definition}
Now, $(k^*)^{\Delta(1)}=Hom(\mathbb{Z}^{\Delta(1)},k^*)$ acts diagonally on ${\rm Cox(X)}$ and thus induces an action of the group $G=Hom({\rm Cl}(X),k^{*})$ via the inclusion
\[0\rightarrow Hom({\rm Cl}(X),k^*)\rightarrow Hom(\mathbb{Z}^{\Delta(1)},k^*)\rightarrow Hom(M,k^*)\rightarrow 0\]
obtained from applying $Hom(\_,k^*)$ to the sequence
\[
0\rightarrow M\rightarrow \mathbb{Z}^{\Delta(1)}\rightarrow {\rm Cl}(X)\rightarrow 0.
\]

\begin{teo} The toric variety $X$ is the categorical quotient of the open set $\Sp(Cox(X))\setminus V(\irr_X)$ by the action of $G$.
\end{teo}
\begin{proof} The scheme $\Sp(\Cox 
(X))\setminus V(\irr_X)$ is covered by the affines $\Sp(\Cox 
(X)_{v^{\overline\sigma}})$ as $v^{\overline \sigma}$ runs over the generators of $\irr_X$. The quotient by $G$ of each of these is the ring of invariants (i.e. the degree $0$ part).
If $\sigma$ is any cone and $u\in \sigma^{\vee}\cap M$ then $u$ determines an element  
\[
v_1^{\langle u,v_1\rangle}\dots ~v_d^{\langle u, v_d\rangle}\in\Cox 
(X)_{v^{\overline\sigma}}
\]
because $\langle u,v_i\rangle\geq 0$ if $v_i\in \sigma$. This element has degree $0$ since $[div(x^u)]=0$. This correspondence is an isomorphism between $U_{\sigma}$ and $\Sp\bigl((Cox(X)_{v^{\overline\sigma}})_0\bigr)$ compatible with the glueing maps and thus an isomorphism between $X(\Delta)$ and the categorical quotient.
\end{proof}

\begin{examples}
\ \vspace{2mm}

\begin{itemize}\itemsep 5mm
\item
If $X=\mathbb{P}^n$ the Cox ring is the polynomial ring in $(n+1)$ generators, one for each ray of its fan and with $deg(v_i)=1$, for all $i$. Since every subset of $n$ rays spans a cone, the irrelevant ideal is $J_X=(v_0,\dots, v_n)$. We thus recover the usual construction of $\mathbb{P}^n$ as a quotient of $\mathbb{A}^{n+1}\setminus\{0\}=\Sp(k[v_1,\dots, v_d])\setminus V(\irr_X)$ by the action of $k^{*}$.

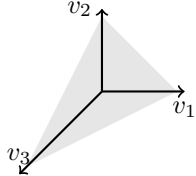
\begin{center}
\begin{tikzpicture}[scale=1]
	\shade[top color=gray!20, bottom color=gray!20] (0,0) -- (1,0) -- (0,1);
	\shade[top color=gray!20, bottom color=gray!20] (0,0) -- (0,1) -- (-1,-1);
	\shade[top color=gray!20, bottom color=gray!20] (0,0) -- (-1,-1) -- (1,0);
	\draw[->,thick] (0,0) -- (1.1,0) node[below] {$v_1$};
	\draw[->,thick] (0,0) -- (0,1.1) node[left] {$v_2$};
	\draw[->,thick] (0,0) -- (-1.1,-1.1) node[above] {$v_3$};
\end{tikzpicture}
\end{center}

\item
Now, for a singular example consider the fan with a single cone generated by $v_1=2e_1-e_2$ and $v_2=e_2$. The resulting toric variety is singular because the determinant of $v_1, v_2$ is not $\pm 1$.
Its Cox ring is $k[v_1,v_2]$ with a $\mathbb{Z}/2\mathbb{Z}$-grading given by $deg(v_j)=1$. The irrelevant ideal is $\irr_X=(1)$ and $X=\Sp(Cox(X)^{G})$ which is the subalgebra of the polynomial ring generated by the even degree pieces. It is thus isomorphic to $k[x^2,xy,y^2]$ that is, the cone over a smooth quadric.

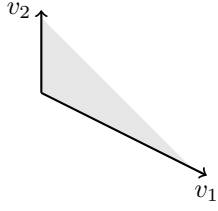
\begin{center}
\begin{tikzpicture}[scale=1]
	\shade[top color=gray!20, bottom color=gray!20] (0,0) -- (2,-1) -- (0,1);
	\draw[->,thick] (0,0) -- (2.2,-1.1) node[below] {$v_1$};
	\draw[->,thick] (0,0) -- (0,1.1) node[left] {$v_2$};
\end{tikzpicture}
\end{center}
More generally, when the fan $\Delta$ has only one cone the above construction shows that any affine toric variety is a quotient of some affine space by a reductive group. In particular, every affine toric variety is Cohen-Macaulay.

\item
Consider the fan $\Delta_n$. The corresponding toric variety is the Hirzebruch surface $\mathbb{F}_n$. Its Cox ring is the $\mathbb{Z}^2$-graded polynomial ring $k[v_1,\dots, v_4]$ with $deg(v_1)=deg(v_3)=(1,0)$, $deg(v_2)=(-n,1)$ and $deg(v_4)=(0,1)$. The irrelevant ideal is $J_X=(v_1v_2, v_2v_3, v_3v_4, v_4v_1)$.
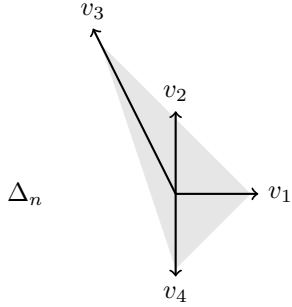
\begin{center}

\begin{tikzpicture}[scale=1]
	\shade[top color=gray!20, bottom color=gray!20] (0,0) -- (0,1) -- (1,0);
	\shade[top color=gray!20, bottom color=gray!20] (0,0) -- (0,1) -- (-1,2);
	\shade[top color=gray!20, bottom color=gray!20] (0,0) -- (1,0) -- (0,-1);
	\shade[top color=gray!20, bottom color=gray!20] (0,0) -- (-1,2) -- (0,-1);
      
	\draw[->,thick] (0,0) -- (1.1,0) node[right] {$v_1$};
	\draw[->,thick] (0,0) -- (0,1.1) node[above] {$v_2$};
	\draw[->,thick] (0,0) -- (0,-1.1) node[below] {$v_4$};
	\draw[->,thick] (0,0) -- (-1.1,2.2) node[above] {$v_3$};
	
	\draw (-2,0) node {$\Delta_n$};
\end{tikzpicture}

\end{center}

\end{itemize}
\end{examples}

The following theorem motivates the general definition of ${\rm Cox(X)}$ (see ~\cite{hk})
\begin{teo} \label{motiv} For any $T$-invariant Weil divisor $D$, the vector space of global sections of the sheaf associated to $D$,  $H^0(X,D)$ is isomorphic to the degree $[D]$ part of $\Cox(X)$. 
\end{teo}
\begin{proof} Suppose $D=a_1Y(v_1)+\dots+a_kY(v_k)$. A rational function $x^u$ is in $H^0(X,D)$ iff $div(x^u)+D\geq 0$ or equivalently if $\langle u,v_i\rangle \geq -a_i$ for all $i$. As a result, the expression $v_1^{\langle u,v_1\rangle+a_i}\dots, v_d^{\langle u,v_d\rangle+a_d}$ has only nonnegative exponents and thus is a monomial of degree $[D]$ in ${\rm Cox(X)}$. This correspondence is obviously bijective. Since the sections of the form $x^u\in H^0(X,D)$ are a basis, the statement follows.\end{proof}
When $X$ is a projective toric variety, the identification between graded components of $\Cox 
(X)$ and global sections of line bundles allows us to give a natural geometric interpretation of the ideal $\irr_X$:
\begin{lemma} Let $D$ be any ample line bundle on $X$, then 
\[
\irr_X = \sqrt{\big(H^0(X,D)\big)}
\]
where the right side of the equality is the ideal generated by the images of the sections of $H^0(X,D)$ in ${\rm Cox(X)_D}$ under the above correspondence. 
\end{lemma}
\begin{proof} Since $D=a_1Y(v_1)+\dots+a_dY(v_d)$ is ample, it follows that for every cone $\sigma$ there are sections $x^{u(\sigma)}\in H^0(X,D)$ such that $\langle u,v_i\rangle\geq -a_i$ for all $i$ with equality if and only if $v_i\in \sigma$. Thus the monomial $v_1^{\langle u(\sigma),v_1\rangle+a_1}\dots v_d^{\langle u(\sigma),v_d\rangle+a_d}$ involves precisely the rays which do not belong to $\sigma$. Thus its radical equals $v^{\overline{\sigma}}$ and the two ideals above coincide.
\end{proof}

\begin{example} The Hirzebruch surface $\mathbb{F}_n$ has the ample line bundle $D=Y(v_2)+(n+1)Y(v_3)$. Its global sections are the integral points in the convex hull of $1,x^{1+n},y^{-1}$ and $xy^{-1}$. These sections correspond to the monomials $v_2v_3^{1+n}$, $v_1^{1+n}v_2$, $v_4v_3$ and $v_1v_4$ in $\Cox 
(\mathbb{F}_n)_{D}$. The radical of the ideal they generate equals $J_X$.\\
\end{example}

\section{Algebraic varieties and their Cox rings}

In this section we introduce the total coordinate ring or Cox ring of an algebraic variety $X$, defined over an algebraically closed field $k$. We will assume that our variety $X$ satisfies:
\begin{itemize}
\item{The group ${\rm Pic}(X)$ is free of rank $r$ and coincides with ${\rm Cl}(X)$}  
\end{itemize}
Note that if $X$ is smooth and $\Pic(X)$ is free the equality is satisfied.

\begin{remark}
It is possible to define the Cox ring also when $\Pic(X)$ differs from $\Cl(X)$: see~\cite{bh2,COX,EKW} for a definition done via $\Cl(X)$ and~\cite{bh,ka} for a definition based on $\Pic(X)$. In the first case the Cox ring is factorial~\cite[Prop. 8.4]{bh}.
\end{remark}

\subsection{Cox rings}\label{sec: Coxrings}

Let $L_1, \dots, L_n$ be a collection of line bundles whose classes form a basis for ${\rm Pic(X)}$. 
We define the {\bf\em Cox ring} of $X$ to be:
\[
\Cox(X,L) := \bigoplus_{(m_1,\dots, m_n)\in\zz^n} H^0(X,m_1L_1+\dots+m_nL_n).
\]
with the multiplication induced by the product of rational functions.

Different choices of bases for $\Cl(X)$ yield (non-canonically) isomorphic Cox rings and we will adopt the notation $\Cox(X)$ whenever it is not important to stress the dependence on the basis $L$.

The Cox ring of $X$ is naturally graded by the Picard group of $X$. This means that there is a $\zz$-module map:
\[
\deg: \Cox(X)\map \Cl(X)
\]
which sends a section $x^D$ to its class $[D]$. Note that the image of the map ${\rm deg}$ is precisely the set $\Eff(X)$ of classes of effective divisors.

\subsection{The torus action}

Denote by $T_X$ the algebraic torus ${\rm Hom}({\rm Pic(X)},k^*)\cong{k^*}^r$.  This torus acts naturally on the Cox ring of $X$. If $x$ is a section in $H^0(X,D)$, where $D = a_1 L_1+\dots + a_r L_r$, then the action is given by
\[
(t_1,\dots,t_r)\cdot x := t_1^{a_1}\cdots t_r^{a_r}~x.
\]
The action extends to an action of $T_X$ on the affine scheme $\overline{X}=\Sp(Cox(X))$.

\begin{example} By Theorem~\hyperref[motiv]{\ref{motiv}} the Cox ring defined above coincides with Cox's construction in the case of toric varieties which satisfy the above restrictions. Moreover, by Corollary 2.10 in~\cite{hk}  the ring ${\rm Cox(X)}$ is a polynomial ring only if $X$ is a toric variety.
\end{example}
\begin{example}
If $X$ is a projective variety and its homogeneous coordinate ring $R$ in some embedding is a unique factorization domain then $\Cl(X)=\mathbb{Z}$ and $\Cox 
(X)=R$. In particular we see that the Cox rings of all grassmannians coincide with their homogeneous coordinate rings in the Pl\"ucker embedding. 
\end{example}

We will discuss the Cox rings of other homogeneous spaces in Section~\hyperref[homogeneous]{\ref{homogeneous}}.

\begin{example}
\label{ex:DP5}
Let $S$ be the smooth surface obtained by blowing up $\pp^2$ at $p_1,\ldots,p_4$ distinct points, no three of which are collinear. 
In this case $\Pic(S)$ is a free abelian group of rank $5$ generated by the pullback $H$ of a line in $\mathbb{P}^2$ and by the four exceptional divisors $E_1,\dots, E_4$ (see~\cite{bp}).

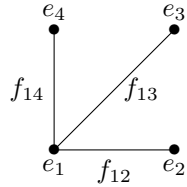
\begin{center}
\begin{tikzpicture}[scale=.8]
    \draw (0,0)node{$\bullet$}node[below] {$e_1$} --node[below]{$f_{12}$} (2,0);
    \draw (2,0)node{$\bullet$}node[below] {$e_2$};
	 \draw (2,2)node{$\bullet$}node[above] {$e_3$};
    \draw (0,0) -- node[left]{$f_{14}$}(0,2);
    \draw (0,2)node{$\bullet$}node[above] {$e_4$};
    \draw (0,0) -- node[right]{$f_{13}$} (2,2);
\end{tikzpicture}
\end{center}

$\Cox(S)$ has $10$ generators corresponding to the $(-1)$-curves on $S$
which we denote by $e_1,\dots,e_4,f_{12},\dots,f_{34}$, graded by $\deg(f_{ij})=H-E_i-E_j$, $\deg(e_i)=E_i$. The ideal of relations is generated by the five quadrics
\[
\begin{array}{cc}
f_{14}f_{23}-f_{12}f_{34} -f_{13}f_{24} & e_2f_{12} -e_3f_{13}-e_4f_{14},\\
e_1f_{12} - e_3f_{23}-e_4f_{24} & e_1f_{13} -e_2f_{23}+e_4f_{34},\\
e_1f_{14} -e_2f_{24}-e_3f_{34}
\end{array}
\]

Thus ${\rm Cox(S)}$ is isomorphic to the homogeneous coordinate ring of the grassmannian ${\rm Gr}(2,5)$ of lines in $\mathbb{P}^4$ in its Pl\"ucker embedding (note however that this is not an isomorphism of graded rings). 
\end{example}

The variety $S$ is an instance of a Del Pezzo surface. We will describe their Cox rings in Section~\hyperref[Del Pezzo]{\ref{Del Pezzo}}. Moreover, Del Pezzo surfaces are instances of Fano varieties  (i.e. those varieties whose anticanonical divisor $-K_X$ is ample). By Corollary 1.3.1 in~\cite{bchm} the Cox rings of all Fano varieties are finitely generated. It is an interesting open problem to describe generators and relations for their Cox rings.

Note that both the polynomial ring and the homogeneous coordinate rings of Grassmannians in their Pl\"ucker embedding are unique factorization domains. This property holds for the Cox rings of all normal projective varieties with finitely generated Cox rings by~\cite{EKW}.

\begin{propo}\label{gen}
Let $D$ be an integral Cartier divisor of $X$ such that $ h^0(X,D) = 1$ then $x^D$ is a generator of $\Cox(X)$.
\end{propo}
\begin{proof}
Write $x^D$ as a polynomial $P$ in the generators $e_1,\dots,e_r$ of $\Cox(X)$ and let $m = e_1^{a_1}\dots e_r^{a_r}$ be a monomial of $P$.  This means that $D$ is linearly equivalent to $a_1E_1+\dots+a_rE_r$
where $E_i$ is the divisor of equation $e_i=0$. By hypothesis this can happen only if $D = E_i$ for some $i$.
\end{proof}

Observe that if $C$ is an integral curve on an algebraic surface $S$ such that $C^2<0$, then $h^0(S,C) = 1$ (if $C'\in |C|$ is distinct from $C$ then $C^2=C'\cdot C\geq 0$).
By Proposition~\hyperref[gen]{\ref{gen}} we get that $x^C$ appears in any set of generators of $\Cox(S)$.

\begin{example}
Let $X$ be the blow-up of $\pp^2$ at nine points in the intersection of two general plane cubics then $\Cox(X)$ is not finitely generated because $X$ contains infinitely many $(-1)$-curves. The same happens if $X$ is a K3 surface with $\rk\Pic(X) = 20$, since in this case there are infinitely many $(-2)$-curves on $X$.

An obvious necessary condition is that the cone of effective divisors must be polyhedral. The condition is not sufficient since there are varieties with finitely generated effective cones whose Cox rings are not finitely generated (see Example 1.8 in~\cite{ht}).
\end{example}

\subsection{Cox rings and GIT}
In what follows we will assume the Cox ring of $X$ to be finitely generated.

Let $\overline{X}:=\Sp(\Cox(X))$; recall that the torus $T_X={\rm Hom}({\rm Pic}(X),k^*)\cong (k^*)^r$ has character lattice ${\rm Pic}(X)=\mathbb{Z}^r$ and acts on $\overline{X}$.
Let $D$ be a divisor on $X$ and
\[
R_D := \bigoplus_{m=0}^\infty\Cox(X)_{mD}.
\]
The inclusion $R_D\subseteq\Cox(X)$ induces a rational map:
\[
\pi_D: \overline{X}\rmap \Proj(R_D)
\]
which is constant on $T_X$-orbits.
It is thus natural to expect that quotients of $\overline{X}$ by the action of $T_X$ may yield some insight to the geometry of $X$.

In this section we will survey results from~\cite{hk} describing the relationship between $X$ and quotients of $\overline{X}$ by $T_X$. We begin by showing that, in a manner entirely analogous to the case of toric varieties, $X$ can be obtained as a quotient of an open subset of $\overline{X}$ by $T_X$. First we recall some basic definitions from Geometric Invariant Theory.

The linearization of the trivial line bundle ${\mathscr L}$ on $\overline{X}$ corresponding to the character $D$ is the morphism 
\[\Cox(X)[t]\rightarrow \Cox(X)[t]\otimes_kk[\mathbb{Z}^r]\] 
which extends the action of $T_X$ on $\Cox(X)$ and maps $t$ to $t\otimes x^{D}$.
We use $\mathscr{L}_{D}$ to denote the trivial line bundle together with the linearization specified by $D$. A $T_X$-invariant section $s$ of $\mathscr{L}_D$ is a $\Cox(X)$- algebra homomorphism $s:\Cox(X)[t]\rightarrow \Cox(X)$ which makes the following diagram commute.
\begin{center}
\begin{tikzpicture}[scale=.8]
\draw (0,0) node (a){$\Cox(X)[t]$};
\draw (4,0) node (b){$\Cox(X)[t]\otimes k[\mathbb{Z}^r]$};
\draw (0,-2) node (c){$\Cox(X)$};
\draw (4,-2) node (d){$\Cox(X)\otimes k[\mathbb{Z}^r]$};
\draw [->] (a)--(b);
\draw [->](a)--node[left]{$s$} (c);
\draw [->](b)--node[right]{$s\times id$} (d);
\draw [->] (c)--(d);
\end{tikzpicture}
\end{center}
Thus, we have a bijection between the $T_X$-invariant sections of $\mathscr{L}_D$ on $\overline{X}$ and the $D$-th graded component of the Cox ring, 
\[
H^0(\overline{X},\mathscr{L}_D)^{T_X}=\Cox(X)_D.
\]
The higher tensor powers of $\mathscr{L}$ inherit a natural linearization with $\mathscr{L}_D^{\otimes n}=\mathscr{L}_{nD}$. As a result, the invariants of the ring of sections of $\mathscr{L}$ with respect to the linearization $D$ can be identified with a subalgebra of $\Cox(X)$
\[
\bigl(\bigoplus_{m=0}^\infty H^0(\overline{X},\mathscr{L}_D^n)\bigr)^{T_X}=R_D.
\]
The  GIT quotient of $\overline{X}$ by $T_X$ (via the linearization determined by $D$) is the scheme $\Proj(R_D)$. A point $x\in\overline{X}$ is semistable if there exists a $T_X$-invariant section $s$ of some tensor power of $\mathscr{L}_D$  such that $s(x)\neq 0$. The scheme $\Proj(R_D)$ is a good quotient of the set of semistable points in $\overline{X}$ by $T_X$ (note that the set of semistable points and hence the quotient depends on the linearization). 

We now introduce the analogue of the irrelevant ideal of toric varieties.
Given a line bundle $D$ in $X$ consider the ideal of $\Cox(X)$ defined by
\[
\irr_{X,D} := \sqrt{B_{X,D}}\hspace{5mm}\text{where}\hspace{5mm} B_{X,D} := (R_D).
\]
Note that $V(\irr_L)$ consists of those points in $\overline{X}$ which are not semistable with respect to the linearization determined by $L$.

The {\bf\em irrelevant ideal} of $X$ is the ideal $\irr_X:=\irr_{X,L}$ for some ample line bundle $L$ on $X$. Note that $\irr_X$ is independent of the choice of ample line bundle.

\begin{propo}\label{quotient} The variety $X$ is a good geometric quotient of $\hat{X} := \Sp(\Cox(X))\setminus V(\irr_X)$ by $T_X$.
\end{propo}
\begin{proof} Let $L$ be a very ample line bundle on $X$. The GIT quotient of $\Cox(X)$ by $\mathscr{L}_L$ is 
\[
\Proj(\bigoplus_{m=0}^\infty \Cox(X)_{mL})=\Proj(R_L)\cong X,
\]
where the first equality follows since the $A$-th graded component of the Cox ring is the space of global sections $H^0(X,A)$. Since the semistable points of $\overline{X}$ are $\hat{X}$ the statement follows. Moreover, as shown in~\cite[Proposition 2.9]{hk}, this quotient is a geometric quotient.
\end{proof}

\begin{remark}
Proposition~\hyperref[quotient]{\ref{quotient}} shows that $X$
is completely determined by its Cox ring and some combinatorial data in the grading group of $\Cox(X)$. The combinatorial data fixes the ``small birational type" of $X$.
In the case of a projective $X$ one may take its GIT chamber as combinatorial data. To include also nonprojective $X$, one needs more generally the``bunches of cones" of~\cite{bh2}.
\end{remark}

\begin{example}
For the Del Pezzo surface from Example~\hyperref[ex:DP5]{\ref{ex:DP5}}, we have that 
\[
\Cox(S)_{-K_S}=H^0(S,-K_S)= \langle e_1e_2f_{12}^2f_{34},\dots,e_3e_4f_{34}^2f_{12}\rangle.
\]
Moreover, $-K_S$ is very ample so every section of the higher tensor powers $-mK_S$ is a product of sections in $-K_S$. As a result, the irrelevant ideal of $S$ is the radical of the ideal of $\Cox(S)$ generated by the above sections.
\end{example}

In view of proposition~\hyperref[quotient]{\ref{quotient}} it is natural to ask the following question: What are the possible quotients of $\overline{X}$ as we vary the linearization $D$ and how do these spaces relate to $X$? The answer is that the cone of effective divisors on $X$ admits a decomposition into finitely many polyhedral chambers in whose interior the GIT quotient is constant. Moreover, there are rational maps from $X$ to each of these quotients.

A {\bf\em chamber} $\cC$ of $X$ is the set of all the effective divisors $D$ having the same $J_D$. (or equivalently the set of divisors whose linearizations induce a fixed GIT quotient).

It follows from the definition that $J_L = J_{nL}$ and that if $L, M$ are two line bundles such that $J_L = J_M$ then $J_{aL+bM} = J_L$ for any pair of positive integers $a,b$. 

The decomposition into chambers is a coarsening of one induced from toric geometry as follows. Suppose $k[A]/I_{X}=\Cox(X)$ is a presentation for the Cox ring of $X$ and that $D$ is an ample divisor on $X$. Then $X$ is a subscheme of the toric variety $Y:=\Proj(\oplus_{m\in \mathbb{N}}k[A]_{mD})$ and we will show that their birational geometry is closely related. Note that $X$ and $Y$ have the same effective cone and that moreover:
\begin{lemma} Any chamber of $X$ is a union of chambers of $Y$.
\end{lemma}
\begin{proof} For any divisor $D\subset X$ consider the ideals of $k[A]$ and $\Cox(X)$ generated by the $[D]$-graded part of the two rings:
\[
B_{Y,D}:=(\bigoplus_{m\in\mathbb{N}} k[A]_{mD}),\hspace{1cm} B_{X,D}:=(\bigoplus_{m\in\mathbb{N}} Cox(X)_{mD}).
\]
If $\sqrt{B_{Y,D}}=\sqrt{B_{Y,F}}$ in ${\rm k[A]}$ then  
\[\sqrt{B_{Y,D}+I_X}=\sqrt{\sqrt{B_{Y,D}}+I_X}=\sqrt{\sqrt{B_{Y,F}}+I_X}=\sqrt{B_{Y,F}+I_X}\]
so $J_D=\sqrt{B_{X,D}}=\sqrt{B_{X,F}}=J_F$ in $\Cox 
(X)$ and the statement follows. 
\end{proof}
Note that the chamber decomposition of $Y$ is finite (the ideal $\sqrt{B_{Y,D}}$ is monomial and there are finitely many radical monomial ideals in any polynomial ring) so it follows that the chamber decomposition of $X$ is also finite:
\[
\Eff(X) := \cC_1\cup\dots\cup\cC_l
\]
The chamber decomposition for toric varieties $Y$ has been studied. For a purely combinatorial description see~\cite{op}. The previous proposition provides an algorithm for determining the chambers once an explicit presentation for the Cox ring of $X$ is known.

\begin{example}
Let $S$ be as in Exercise~\hyperref[ex:DP5]{\ref{ex:DP5}}, i.e. the blow-up of $\pp^2$ at four non-collinear points, and let $D = 11H-5E_1-3E_2-2E_3-E_4$. Observe that $D$ is ample since it has positive intersection with all the $(-1)$-curves of $S$.
The radical of the monomial ideal generated by a basis of $H^0(S,D)$ can be calculated by means of the following Macaulay 2 script:
{\small \begin{verbatim}
i1 : R=QQ[x,y,z,w,a,b,c,d,e,f,Degrees=>{{0,1,0,0,0},{0,0,1,0,0},
{0,0,0,1,0},{0,0,0,0,1},{1,-1,-1,0,0},{1,-1,0,-1,0},{1,-1,0,0,-1},
{1,0,-1,-1,0},{1,0,-1,0,-1},{1,0,0,-1,-1}},Heft=>{3,1,1,1,1}]

i2 : radical monomialIdeal basis({11,-5,-3,-2,-1},R)
\end{verbatim} }
The calculation shows that each monomial contains $5$ of the $10$ variables. This implies that $Y:=\Proj(\oplus_{m\in \mathbb{N}}k[A]_{mD})$ is a $\qq$-factorial toric variety because all the cones of its defining fan are simplicial.
Observe that if instead of $D$ one choses $-K_S$, then the corresponding $Y$ would no longer be $\qq$-factorial.

Observe that the toric variety $Y$ is not projective since all of its $1$-rays live in the half-space $3x_1+x_2+x_3+x_4+x_5\geq 0$ so that its defining fan is not complete. Since the $1$-rays are the same for all the toric varieties with the same Cox ring as $Y$, then none of them is projective.
\end{example}

\begin{example}
The effective cone of an Hirzebruch surface $\ff_n$ with $n\geq 2$ is generated by the $(-n)$-curve, whose class is $[v_2]$ and by a fiber, whose class is $[v_1]$, of the rational fibration. An effective divisor $D$ of class $a[v_1]+b[v_2]$ is ample iff $a-nb > 0,\ b > 0$. 
\begin{center}
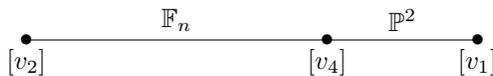
\begin{figure}[h]
\begin{tikzpicture}[scale=2]
	\draw[-] (-2,0) node {$\bullet$} node[below] {$[v_2]$} --node[above] {$\ff_n$} (0,0) node {$\bullet$} node[below] {$[v_4]$} -- node[above] {$\pp^2$} (1,0) node {$\bullet$} node[below] {$[v_1]$} ;
\end{tikzpicture}
\caption{Mori chambers of $\ff_n$ with the corresponding values of $\Proj(R_D)$.}
\end{figure}
\end{center}
\end{example}

The chambers can be interpreted in terms of the geometry of $X$ (see ~\cite{hk}). 
If $D$ is an effective divisor there is an induced rational map:
\[
f_D:X\rmap \Proj(R_D).
\]
We say that two divisors $D_1$ and $D_2$ are {\bf\em Mori equivalent} if there is an isomorphism between their images which makes the natural triangular diagram commutative. If $\Cox(X)$ is finitely generated then the GIT-Chambers and the Mori chambers of $X$ coincide. In particular, the Mori chambers of $X$ are a coarsening of the Mori chambers of a toric variety $Y$, as in the previous proposition. A simple direct construction for
the the chamber decomposition of a torus action on a factorial affine variety is given in~\cite[Section 2]{bh0}.

\section{Mori Dream Spaces}
In this section we give a brief description of the Mori program and state its relationship with Mori Dream Spaces. Given a $\mathbb{Q}$-factorial algebraic variety $X$, the minimal model program is a procedure to modify $X$ within its birationality class with the aim of making its canonical divisor $K_X$ nef (i.e. ensuring that $K_X\cdot c\geq 0$ for every curve $c$).
 
At each step, curves on which $K_X$ is negative are contracted as follows: if $K_X$ is not nef, the cone theorem (see~\cite{D} pag. 144) shows that there exists an (extremal) rational curve $z$ with $K_X\cdot z<0$. By the contraction theorem (see~\cite{D} pag.153) there exists a morphism with connected fibers $f:X\rightarrow Y$ which contracts $z$. Moreover, one of the following occurs:
\begin{enumerate}\itemsep 2mm
\item{The morphism $f$ is birational and there is a divisor $C$ in $X$ such that ${\rm \dim(f(C))<\dim(C)}$. In this case $Y$ is $\mathbb{Q}$-factorial and $f$ is called a {\bf\em divisorial contraction}.}
\item{The morphism $f$ is a small contraction (i.e. an isomorphism in codimension one). In this case $Y$ is no longer $\mathbb{Q}$-factorial and (conjecturally if ${\rm dim(X)}>3$) there exists a unique {\bf\em flip} of $f$, that is, a $\mathbb{Q}$-factorial variety $X^{+}$ and a different small contraction $f^{+}:X^{+}\rightarrow Y$ such that $K_{X^{+}}$ is $f^{+}$-ample (i.e. ${K_{X^+}}_{|f^{-1}(U)}$ is ample for every open affine $U\subseteq Y$).
Note that in this case the composition of ${f^{+}}^{-1}\circ f:X\rightarrow X^{+}$ is birational.}
\item{$\dim(Y)<\dim(X)$ and $f:X\rightarrow Y$ is called a {\bf\em Mori fiber space}. If $X$ is smooth, the fibers of this map are Fano varieties.}
\end{enumerate}
We can thus construct sequences of birational maps
\[X=X_0\rightarrow X_1\rightarrow X_2\rightarrow \dots \rightarrow X_m\]
where the $X_i$ are $\mathbb{Q}$-factorial varieties and $g_i:X_{i}\rightarrow X_{i+1}$ is either a a divisorial contraction or a composition ${f^{+}}^{-1}\circ f$ coming from a flip.
 
Moreover, {\it under the additional assumption that there are no infinite sequences of flips}, every sequence as above must eventually stop. Thus, after finitely many steps we reach a variety $X_m$ which is either a Mori fiber space (and we have essentially reduced the study of the geometry of $X_m$ to that of fibrations above lower dimensional varieties) or has nef canonical divisor. In this case $X_m$ is a {\bf\em Minimal Model}.

The questions of existence and termination of flips are the main obstacles for the minimal model program in dimension greater than $3$ (for dimension $3$ and characteristic $0$ they are theorems due to Mori). These questions can be answered positively in the case of varieties $X$ with finitely generated Cox rings. They are called {\bf\em Mori Dream Spaces} because of the following theorem (see ~\cite{hk})
\begin{teo} If $\Cox 
(X)$ is finitely generated then the Mori program can be carried out for any divisor on $X$. The required flips and contractions exist and every sequence terminates.
\end{teo}

In the next section we will see how to construct a $D$-flip by moving from one Mori chamber to the other of the effective cone of $X$.

\subsection{Small birational maps}

We conclude our description of Cox rings by studying how they relate with {\bf\em flips} and {\bf\em small birational maps} in general. A birational map $f: X\rmap Y$ is {\em small} if there are two open subsets $U_X\subseteq X$ and $U_Y\subseteq Y$, whose complements have codimension $\geq 2$ and such that $f_{|U_X}: U_X \map U_Y$ is an isomorphism. 
\begin{propo}\label{small}
If $f: X\rmap Y$ is a small birational map, then $\Cox(X)\cong\Cox(Y)$.
\end{propo}

Propositions~\hyperref[small]{\ref{small}} and~\hyperref[quotient]{\ref{quotient}} imply that {\bf any} such $Y$ is the image of a $\pi_D$ for some $D$ in a $\cC_i\subset\Eff(X)$.
In particular this implies that $X$ has finitely many flip images.

\begin{example}\label{flip}

Let $\pi: X\map\pp^3$ be the blow-up of $\pp^3$ at two distinct points, then there is a smooth toric variety $Y$ and a small birational map
\[
f: X\rmap Y
\]
with $Y$ non-isomorphic to $X$. Consider a basis $v_1, v_2, v_3$ of $\zz^3$ and the vectors $v_4 = -(v_1+v_2+v_3),\ v_5 = -v_3,\ v_6 = -v_2$. The fans of $X$ and $Y$ are:
\[
\begin{array}{c|c}
X & Y \\
\hline
(1,2,3) & (1,2,3) \\
(1,2,5) & (1,2,5) \\
{\bf (1,4,5)} & {\bf (1,5,6)} \\
(2,4,5) & (2,4,5) \\
(1,3,6) & (1,3,6) \\
{\bf (1,4,6)} & {\bf (4,5,6)} \\
(3,4,6) & (3,4,6) \\
(2,3,4) & (2,3,4) 
\end{array}
\]
where the notation $(1,2,3)$ means the three dimensional cone generated by $v_1, v_2, v_3$. Observe that since
\[
v_1+v_4 = v_5 + v_6
\]
there are two ways of triangulating the quadrilateral $(1,4,5,6)$. The varieties $X$ and $Y$ correspond to each of the two triangulations. Since the birational map between the two is small, it follows that $\Cox(X)\cong\Cox(Y)$. 

The cone of effective divisors of $X$, is generated by the exceptional divisors $D_5, D_6$ and the strict transform $D_1$ of a hyperplane in $\pp^3$ passing through the two points. The {\em cone of curves} of $X$ is generated by the classes $e_{15}, e_{16}$ of lines in $D_5, D_6$ and the class of the strict transform of the line $e_{14}$ through the two points. Looking at the intersection between curves and divisors
\[
\begin{array}{c|rrr}
	& D_1 & D_5 & D_6 \\
	\hline
e_{15}	& 1 & -1 & 0 \\
e_{16}	& 1 & 0 & -1 \\
e_{14}	 & -1 & 1  & 1 
\end{array}
\]
we obtain that $aD_1+bD_5+cD_6$ is ample if and only if $\max(b,c) < a < b+c$. This region is represented by the shaded triangle of the following picture.
\begin{center}
\begin{figure}[h]\label{fig-chambers}
\begin{tabular}{cc}
\begin{tikzpicture}[scale=1.5]
	\shade[top color=gray!20,bottom color=gray!20] 
      (1,0) -- (0,1) -- (.66,.66);
	\draw (0,0) node {$\bullet$};
	\draw[-] (0,0) --  (0,2);
	\draw[-] (0,0) -- (2,0);
	\draw[-] (2,0) -- (0,2);
	\draw[-] (1,0) -- (0,2);
	\draw[-] (0,1) -- (2,0);
	\draw[-] (1,0) -- (0,1);
	\draw (2,0)node {$\bullet$};
	\draw (0,2)node {$\bullet$};
	\draw (-.2,-.2) node {$D_1$};
	\draw (2.3,-.2) node {$D_5$};
	\draw (-.2,2.3) node {$D_6$};
\end{tikzpicture}
&
\hspace{1cm}
\begin{tikzpicture}[scale=3.5]
	\shade[top color=gray!20,bottom color=gray!20] 
      (1,0) -- (0,1) -- (.66,.66);
	\draw (0,0) node {$\bullet$};
	\draw[-] (0,0)node[below] {$D_1$} --  (0,1.1);
	\draw[-] (0,0) -- (1.1,0);
	\draw[-, thick] (1,0) -- (0,1);
	\draw[-] (1,0) -- (.5,1);
	\draw[-] (0,1) -- (1,.5);
	\draw[-] (1,0) -- (0,1);
	\draw (.5,.5) node {$\bullet$};
	\draw (.38,.5) node {$Z$};
	\draw (.52,.65) node {$X$};
	\draw[dashed,-] (.63,.63) node {$\bullet$} -- (.35,.35) node {$\bullet$} node[left] {$Y$};
\end{tikzpicture}
\end{tabular}
\caption{Walking through the Mori chambers of $X$. }
\end{figure}
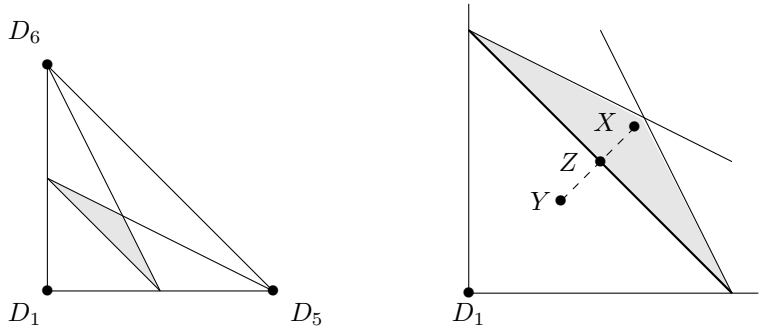
\end{center}

Consider the divisor $D_X = 3D_1+2D_5+2D_6$. Note that $D_X$ is ample on $X$ and denote its coordinate in the effective cone of $X$ by a bullet. Imagine it moves towards the boundary of the region as depicted in the figure. The wall is crossed at $D_Z = D_X+D_Y$ where $D_Y = 3D_1+D_5+D_6$ is not ample on $X$ because
\[
D_Y\cdot e_{14} = -1.
\]
Since $D_Z\cdot e_{14} = 0$ the curve $e_{14}$ is contracted in the {\bf\em categorical} quotient $\hat{X}\map Z = \Proj (R_{D_Z})$. Moving to the new chamber we obtain the  {\bf\em geometric} quotient $\hat{X}\map Y = \Proj (R_{D_Y})$. If we denote by $\phi:=\phi_{|D_Z|}$, we have just witnessed the $D_Y$-flip of $\phi$.
The effect of walking between chambers can be visualized by considering the dual polyhedra of the three toric varieties.

\begin{center}
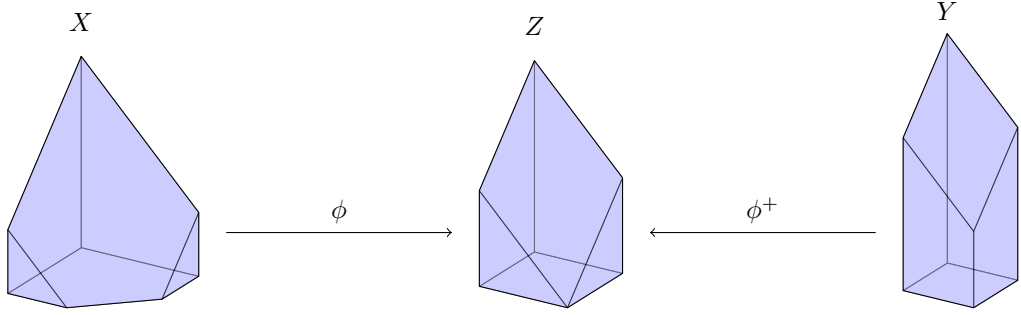
\begin{figure}[h]
\begin{tabular}{ccccc}
\begin{tikzpicture}[join=round]
\filldraw[fill=blue!20,draw=black](0,0)--(-.976,-.61)--(-.195,-.8)--(1.073,-.686)--(1.561,-.381)--cycle;
\filldraw[fill=blue!20,draw= black](0,0)--(0,2.543)--(-.976,.238)--(-.976,-.61)--cycle;
\filldraw[fill=blue!20,draw= black](0,0)--(1.561,-.381)--(1.561,.467)--(0,2.543)--cycle;
\filldraw[fill=blue!20,draw= black,semitransparent](1.073,-.686)--(1.561,-.381)--(1.561,.467)--cycle;
\filldraw[fill=blue!20,draw= black,semitransparent](0,2.543)--(-.976,.238)--(-.195,-.8)--(1.073,-.686)--(1.561,.467)--cycle;
\filldraw[fill=blue!20,draw= black,semitransparent](-.976,.238)--(-.976,-.61)--(-.195,-.8)--cycle;
\draw (0,3) node {$X$};
\end{tikzpicture}
&
\begin{tikzpicture}
\draw[->] (0,1) --node[above] {$\phi$}  (3,1);
\draw (0,0) {};
\end{tikzpicture}
&
\begin{tikzpicture}[join=round,scale=1.5]
\filldraw[fill=blue!20,draw=black](.781,-.19)--(0,0)--(-.488,-.305)--(.293,-.495)--cycle;
\filldraw[fill=blue!20,draw=black](0,0)--(.781,-.19)--(.781,.657)--(0,1.695)--cycle;
\filldraw[fill=blue!20,draw=black](-.488,-.305)--(0,0)--(0,1.695)--(-.488,.543)--cycle;
\filldraw[fill=blue!20,draw=black,semitransparent](.781,.657)--(.293,-.495)--(.781,-.19)--cycle;
\filldraw[fill=blue!20,draw=black,semitransparent](-.488,.543)--(-.488,-.305)--(.293,-.495)--cycle;
\filldraw[fill=blue!20,draw=black,semitransparent](-.488,.543)--(.293,-.495)--(.781,.657)--(0,1.695)--cycle;
\draw (0,2) node {$Z$};
\end{tikzpicture}
&
\begin{tikzpicture}
\draw[<-] (0,1) --node[above] {$\phi^{+}$}  (3,1);
\draw (0,0) {};
\end{tikzpicture}
&
\begin{tikzpicture}[join=round,scale=1.2]
\filldraw[fill=blue!20,draw=black](0,0)--(.781,-.19)--(.781,1.505)--(0,2.543)--cycle;
\filldraw[fill=blue!20,draw=black](0,0)--(-.488,-.305)--(.293,-.495)--(.781,-.19)--cycle;
\filldraw[fill=blue!20,draw=black](0,0)--(-.488,-.305)--(-.488,1.39)--(0,2.543)--cycle;
\filldraw[fill=blue!20,draw=black,semitransparent](.781,1.505)--(.293,.352)--(.293,-.495)--(.781,-.19)--cycle;
\filldraw[fill=blue!20,draw=black,semitransparent](.293,.352)--(-.488,1.39)--(-.488,-.305)--(.293,-.495)--cycle;
\filldraw[fill=blue!20,draw=black,semitransparent](-.488,1.39)--(.293,.352)--(.781,1.505)--(0,2.543)--cycle;
\draw (0,2.8) node {$Y$};
\end{tikzpicture}
\end{tabular}
\caption{Crossing a wall between maximal adjacent Mori chambers.}
\end{figure}
\end{center}
Observe that $Y$ has a new curve $e_{56}$ which is contracted by $\phi^{+}$ and that $(\phi^{-1}\circ\phi^{+})^*D_Y$ is an ample divisor on $Y$. It is also interesting to observe that both the small contractions $\phi$ and $\phi^{+}$ are given by the canonical divisors of $X$ and $Y$. For example $D_Z = -\frac{3}{2} K_X$.
\end{example}

\section{Examples of Mori Dream Spaces}

\subsection{Homogeneous varieties}\label{homogeneous}
The aim of this section is to describe the Cox ring of a homogeneous variety. 
We recall that spherical varieties and even unirational varieties with a complexity one group action have finitely generated Cox rings, the latter follows from~\cite{k}.
For the more general case of {\em wonderful varieties} see~\cite{br}.

Let $G$ be an algebraic {\bf\em reductive} group and let $P\subset G$ be a {\bf\em parabolic} subgroup of $G$, this is any $P$ which contains a {\bf\em Borel} subgroup of $G$. For us $G = \SL(n,\cc)$ and $B$ will be the subgroup of upper-triangular matrices.
The following is a well-known fact.
\begin{propo}
$X = G/P$ is a projective algebraic variety.
\end{propo}

The action of $B$ on $X$ has an open orbit $X^0\subset X$ whose complementary set:
\[
X\setminus X^0 = D_1\cup\dots\cup D_r
\]
is a union of $B$-invariant divisors which are called the {\bf\em colors} of $X$. If $x^D$ is a defining section for $D$, consider the irreducible $G$-module $V_D\subseteq H^0(X,D)$ containing $x^D$.
\begin{propo}\label{gen-homog}
$\Cox(X)$ is generated by bases of all the $V_{D}$ as $D$ varies in the colors of $X$.
\end{propo}

Recall that the Picard group of $G/P$ is isomorphic to the character group of $P$ by means of the following construction. If $\chi: P\map \cc^*$ is a character consider the equivalence relation on $G\times\cc^*$ defined by: $(g,t)\sim (gb^{-1},\chi(b)t)$,
where $g\in G, b\in P, t\in\cc^*$. The quotient space $L_{\chi}$ is an equivariant line bundle on $X$.

The structure of the Cox ring of a homogeneous variety is well understood.

\begin{propo}
The Cox ring of a homogeneous variety $X=G/P$ is the homogeneous coordinate ring of the quasi-affine variety $G/[P,P]$.
\end{propo}

\begin{example}
The {\em complete flag variety} is the quotient $X := G/B$ where $G = \SL(n,\cc)$ and $B\subset G$ is the subgroup of upper triangular matrices.
Observe that 
\[
B \cong [B,B]\times T,
\]
where $T\subset G$ is the group of diagonal matrices. This implies that $\chi(B) = \chi (T)$ so that $\Pic(X)\cong\zz^{n-1}$.
On the other hand, $T$ acts on $\hat{X} = G/[B,B]$ and the geometric quotient of this action is exactly $X$.

As a concrete example consider the case $n=3$ so that $X$ is isomorphic to the incidence variety:
\[
X_{1,2} := \{(x,l)\in\pp^2\times{\pp^2}^*\ |\ x\in l\}
\]
which comes equipped with the two projection maps $\pi_1, \pi_2$. The generators of $\Pic(X)$ are the pull-backs $H_i := \pi_i^*\Osh_{\pp^2}(1)$ and the colors are bases of the two cohomology groups $H^0(X_{1,2},H_1) = \langle x_0,x_1,x_2\rangle$ and 
$H^0(X_{1,2},H_2) = \langle y_0,y_1,y_2\rangle$. This gives:
\[
\Cox(X)\cong\frac{\cc[x_0,x_1,x_2,y_0,y_1,y_2]}{(x_0y_0+x_1y_1+x_2y_2)}.
\]
$\overline{X}$ is an hypersurface in ${\mathbb A}^3\times {\mathbb A}^3$ and the irrelevant ideal is $J_X = J_{X,H_1}\cap J_{X,H_2}$.
\end{example}

\subsection{Cox rings of Del Pezzo surfaces}\label{Del Pezzo}

A Del Pezzo surface $S$ is a surface whose anticanonical divisor is ample (i.e. a Fano surface). A Del Pezzo surface is isomorphic to either $\mathbb{P}^1\times\mathbb{P}^1$ or to the surface obtained by blowing up $\pp^2$ at $0\leq s\leq 8$ general points. In the second case $\Pic(S)\cong \zz^{s+1}$ and it has a natural basis given by the pullback $L$ of the class of a line in $\mathbb{P}^2$ and the $s$ exceptional divisors $E_1,\dots, E_s$. 

When $s\geq 3$ Del Pezzo surfaces are not toric. As shown in~\cite{bp} they are Mori Dream Spaces. 
More precisely, choosing the basis $(L,E_1,\dots, E_s)$ to define the Cox ring we have,

\begin{propo}
If $4\leq s\leq 7$ then the Cox ring of a Del Pezzo surface $S$ is generated by sections supported on the classes of exceptional curves. If $s=8$, besides the $240$ exceptional curves, we need two additional generators, forming a basis for $H^0(S,-K_S)$. 
\end{propo}
In other words there is a presentation $k[V]/I_S\rightarrow{\rm Cox(S)}$ where $k[V]$ is a polynomial ring with one variable for each exceptional curve in $S$ (and two more if $s=8$). We will now describe the ideal of relations $I_S$. 

First note that, if $D$ is a conic bundle (i.e. a divisor $D\in\Pic(S)$ such that $D^2 = 0,\ D\cdot K_S = -2$), it follows from Riemann-Roch that $h^0(S,D) = 2$ so $D$ defines a morphism $\phi_{|D|}: S\map\pp^1$. Moreover this morphism has $s-1$ reducible fibers which are easily seen to be unions of two exceptional curves. Thus there are $s-1$ quadratic monomials of degree $D$ in ${\rm Cox(X)}$ and hence there are $s-3$ linearly independent linear relations among these. The ideal of relations $I_X$ is generated in degrees $D$ which are natural generalizations of conic bundles as shown in~\cite{LV,SS,STV}.

\begin{propo}
The ideal of relations of the Cox ring of a Del Pezzo surface $S$ is generated by all the relations in nef Picard degrees $D$ with anticanonical degree $-K_X\cdot D=2$. 
\end{propo}
In particular this implies that the Cox rings of Del Pezzo surfaces are quadratic algebras as conjectured in~\cite{bp}. It is shown in~\cite{STV, SX} that moreover these algebras admit quadratic Gr\"obner bases.

\subsection{Blow-ups of projective space}

Let $\pi: X_r\map\pp^n$ be the blow-up of projective space at $r$ points lying on a rational normal curve. Let $E_1,\dots,E_r$ be the exceptional divisors and $H := \pi^*\Osh_{\pp^n}(1)$.
\begin{propo}
The Cox ring of $X_{n+2}$ is generated by
\[
\{x_E\ |\ E\text{ is of the form } E_i\text{ or }\ H-E_{i_1}-\ldots-E_{i_n}\}.
\]
\end{propo}
The divisor $D := H-E_1-\dots-E_{n-1}$ has $h^0(X_{n+2},D) = 3$ and plays the same role of that played by a conic bundle for the case $n=2$.

\begin{propo}
$\Cox(X_{n+2})$ is isomorphic to the total coordinate ring of the grassmannian $G(1,n+2)$.
\end{propo}

In the paper~\cite{ct}, the authors prove that the Cox ring of $X_r$ is finitely generated and they exhibit the generators. 

\section{A technique for determining the relations of ${\rm Cox}(X)$}

Suppose we know a set of generators $\cG := \{E_1,\dots,E_r\}$ of the Cox ring of $X$. This determines a surjective homomorphism $\phi: k[\cG]\rightarrow \Cox(X)$. In this section we review the geometric approach of~\cite{LV} to determine the ideal of relations $I_X={\rm ker }\phi$. 

The polynomial ring $R := k[\cG]$ is ${\rm Pic(X)}$-graded. This grading has a positive coarsening (by the existence of an ample line bundle on $X$) and every  finitely generated $\Pic(X)$-graded $R$-module has a unique minimal $\Pic(X)$-graded free resolution. For the module $\Cox(X) = R/I_X$ this resolution 
is of the form
\[
\dots \rightarrow \bigoplus_{D\in \Pic(X)} R(-D)^{b_{2,D}}\rightarrow \bigoplus_{D\in \Pic(X)} R(-D)^{b_{1,D}} \rightarrow R \rightarrow 0,
\]
where the rightmost nonzero map is given by a row matrix whose entries are a set of minimal generators of the ideal $I_X$. Since the differential of the resolution has degree $0$, it follows that $I_X$ has exactly $b_{1,D}(\Cox(X))$ minimal generators of Picard degree $D$.

Let $\mathbb{K}$ be the Koszul complex on $\cG$. Consider the degree $D$ part of the complex $\Cox(X)\otimes_R \mathbb{K}$. Then
\[
{\rm H}_i \big( (\Cox(X)\otimes_R \mathbb{K})_D \big) = \big( {\rm H}_i(\Cox(X)\otimes_R \mathbb{K}) \big)_D =\big( \Tor^R_i (\Cox(X),k) \big)_D = k^{b_{i,D}(\Cox(X))},
\]
\color{black}
where the last two equalities follow since $\Tor_i^R(A,B)$ is symmetric in $A$ and $B$ and the Koszul complex is the minimal free resolution of $k$ over $R$.  Hence we have the equality
\[
\dim_k \big( {\rm H}_i\big( \bigl(\Cox(X) \otimes_R \mathbb{K}\bigr)_D \big) \big) = b_{i,D}(\Cox(X)).
\]

If we are interested solely in the degrees of generators (and not in the higher syzygies of $I_X$) we may restrict our attention to the following complex:
\begin{equation*}
\bigoplus_{1\leq i<j\leq r} {\rm H}^0 \big( X, D-E_i-E_j \big) \xrightarrow{d_2} \bigoplus_{i=1}^{r} {\rm H}^0 \big( X, D-E_i \big) \xrightarrow{d_1} {\rm H}^0 \big( X, D\big),
\end{equation*}
where $d_2$ sends $\sigma_{ij} \in {\rm H}^0 \big( X,D-E_i-E_j \big)$ to $(0,\dots,0,\sigma_{ij}e_j,0,\dots,0,-\sigma_{ij}e_i,0,\dots,0)$ and $d_1$ sends $\sigma_i \in {\rm H}^0 \big( X, D-E_i \big)$ to $\sigma_ie_i$.
A {\bf\em cycle} is an element of $\ker d_1$; a {\bf\em boundary} is an element of ${\rm im}\, d_2$. The {\bf\em support} of a cycle $\sigma = ( \ldots , \sigma_i , \ldots )$ 
is 
\[
|| \sigma ||= \big\{ E_i :  \sigma_i\neq 0 \big\}
\]
and the size of the support is the cardinality of $|| \sigma ||$, denoted  $| \sigma |$. Note that the ideal $I_X$ has no generators of degree $D$ precisely when every cycle is a boundary. 

One can argue that every cycle is a boundary as follows:
\begin{enumerate}
\item{Describe ways in which a divisor may be removed from the support of a cycle using boundaries, at the cost of possibly introducing new divisors in the support of the cycle.}
\item{Apply the constructions in the previous step to all cycles in a systematic way to reduce their support to at most two elements.}
\item{Conclude, using the fact that every cycle whose support has size two is a boundary. This is an immediate consequence of the fact that {\rm Cox(X)} is a unique factorization domain~\cite{EKW}.}
\end{enumerate}
Step one can be described more precisely by introducing the notion of {\bf\em capture move}. This is a pair $(\cS , C )$, where $\cS \subset \cG$, $C \in \cG \setminus \cS$ 
and the map
\[
\bigoplus_{S \in \cS} {\rm H}^0( X, D - S - C) \otimes {\rm H}^0\big ( X, S)\longrightarrow {\rm H}^0 \big( X,D - C)
\] 
induced by tensor product of sections is surjective. We say that $C$ is the \defi{captured curve}, and that $C$ is \defi{capturable for $D$ by} $\cS$.

If $( \cS , C )$ is a capture move and $\sigma_C \in {\rm H}^0 \big( X, D - C\big)$, 
then we have
\[
\sigma_ C = \sum_{S \in \cS} p_s s, \qquad p_s \in {\rm H}^0 \big( X, D - S - C\big),
\]
and thus we obtain
\[
\sigma_C = \sum_{S \in \cS} p_s c + d_2 \Bigl( \sum_{S \in \cS} \varepsilon _{CS} p_s \Bigr)
\]
where $\varepsilon _{CS} \in \{ \pm 1 \}$.  
Hence, if $\sigma$ is a cycle and $( \cS , C )$ is a capture move, then we can modify 
$\sigma$ by a boundary so that $C \notin ||\sigma ||$. Note, however, that we may have 
added $S$ to $|| \sigma ||$, for all $S \in \cS$, so {\it a priori} the size of the support may 
not have decreased. We need to find and apply capture moves in an organized way to 
ensure we are genuinely decreasing the size of the support of a cycle.

\begin{example}\label{cattura} This is an example of a capture move. Assume that $X$ is the cubic surface. Let $A,B,C$ be exceptional curves on $X$. Assume that $A\cdot B=B\cdot C=1$ and that $A$ and $B$ are disjoint. Let $m_D=\min\{D\cdot V:\text{ V is an exceptional curve in $X$}\}$. If $m_D\geq 1$ then $(\{A,C\},B)$ is a capture move. To prove this fact it suffices to show that $H^1(X,D-A-B-C)=0$ and this is easily done with the Kawamata-Viehweg vanishing theorem since the anticanonical divisor of the cubic surface is sufficiently positive. 
\end{example}

In turn, we can organize the various capture moves systematically, as follows. Assume $D \in \Pic(X)$ and let 
\[
\cM := \big( E_1 , \dots , E_n \big)
\]
be a sequence of elements of $\cG$; define $\cS_i := \cG \setminus \{ E_1 , E_2 , \ldots , E_i \}$ 
for all $i \in \{ 1 , \ldots , n \}$.  We say that $D$ is {\bf\em capturable} (by $\cM$) if 
\begin{enumerate}
\item $( \cS_i , E_i )$ is a capture move for all $i \in \{ 1 , \ldots , n \}$, 
\item $|\cS_n|=2$
\end{enumerate}
Note that if $D$ is capturable then we can use the inductive procedure outlined in the previous definition to conclude that any cycle is a boundary. 
\begin{example} If $X$ is the cubic surface, denote with $e_i$ the exceptional divisors, with $f_{ij}$ the strict transform of the line through the i-th and j-th points and with $g_i$ the unique conic through all points except the $i$-th. If $D$ is a divisor with $m_D\geq 1$ then $D$ is capturable. This is because the sequence $\cM$ of curves in the table satisfies the previous definition with respect to the capture move from Example~\hyperref[cattura]{\ref{cattura}}
\begin{center}
\begin{tabular}{c}
$f_{ij}:(i,j)\neq (1,2), g_1, g_2, e_3, e_4, e_5, e_6,$\\
$g_3, g_4, g_5, g_6, f_{12}$\\
\end{tabular}
\end{center}
with $\cS_{25}=\{e_1,e_2\}$.
\end{example}

This method is used to determine the ideals $I_X$ for Del Pezzo surfaces of degree at least two  in~\cite{LV} and of degree one in~\cite{TVV}.

\bibliographystyle{plain}

\end{document}